\documentclass[12pt]{article}       
\usepackage{graphicx}
\usepackage{setspace}
\usepackage{amsmath}
\usepackage{amssymb}
\usepackage{dcolumn}

\newcommand{\hyper}{{}_2F_1}
\DeclareMathOperator*{\argmin}{arg\,min}

\newcommand{\strat}[2]{\ensuremath{\left<#1,#2\right>}}
\newcommand{\rv}[2]{\ensuremath{X_{\left<{#1,#2}\right>}}}
\newcommand{\payoff}[2]{\ensuremath{\mathbb{V}\left(#1,#2\right)}}

\newcommand{\ii}[4]{\ensuremath{\payoff{\strat{#1}{#2}}{\strat{#3}{#4}}}}

\DeclareMathOperator{\Prob}{Pr}

\begin{document}

\title{A liability allocation game}

\author{Robin K. S. Hankin}

\maketitle
\vfill
{\noindent
  Auckland University of Technology\\
  2-14 Wakefield St\\
  Auckland\\
  New Zealand\\
             Tel.: +64(9)9219999x5106\\
              \tt{robin.hankin@aut.ac.nz}           
}
\newpage

\begin{abstract}
The following problem is considered.  Two players are each required to
allocate a quota of~$n$ counters among~$k$ boxes
labelled~$1,2,\ldots,k$.  At times $t=1,2,3,\ldots$ a random box is
identified; the probability of choosing box~$i$ is~$p_i$.  If a player
has at least one counter in the chosen box, she removes one counter
from it; otherwise she takes no action.  The winner is the first
player to remove all her counters.  The game so described may be
modified so that each player simultaneously, but independently,
identifies a box at random.

This paper analyses this deceptively simple game, which has apparently
not been studied in the literature.  Some analytical and numerical
results are then presented, followed by some challenges for further
work.  
\end{abstract}

\section{Introduction}

Two players are each required to allocate a quota of~$n$ counters
among~$k$ boxes labelled~$1,2,\ldots,k$.  At times~$t=1,2,3,\ldots$ a
random box is identified; the probability of choosing box~$i$ is~$p_i$
and these probabilities are common knowledge.  If a player
has~$\geqslant 1$ counters in the chosen box, she removes one counter
from it; otherwise she takes no action.  The winner is the first
player to remove all her counters, except that if both players remove
their respective last counters simultaneously then the game is a draw.
As such the problem may be considered to be a two-person zero-sum
game.

Further to the above, in which the same box is identified for each
player (the ``common throw'' r\'{e}gime), the game may be modified so
that each player simultaneously, but independently, identifies a box
at random (the ``separate throw'' r\'{e}gime).  In this case the two
players may both use the same set of
probabilities~$p_1,p_2,\dots,p_k$, or each player may use her own
probabilities; however, unless otherwise stated in what follows, the
two sets of probabilities will be identical.

Observe that, in contrast to allocation games such as Colonel
Blotto~\cite{roberson2006} in which the units to be allocated are
beneficial, here the units are detrimental.  The game described here
bears the same relation to resource allocation as chore
division~\cite{peterson2002} does to cake division~\cite{brams1995}:
we call it ``Alice's game''.

The box probabilities will be arranged in a nonincreasing order, so
that~$p_1\geqslant p_2\geqslant\ldots\geqslant p_k$.  The initial
allocation of counters~$n_1, n_2,\ldots,n_k$, with~$\sum n_i=n$, is
termed a {\em strategy}, written~$\left<n_1,n_2,\ldots,n_k\right>$;
strategies thus comprise {\em compositions} in the language of
Hankin~\cite{hankin2006}.  Players may be identified with their
strategies when this does not cause confusion.

Writing~$t_S$ for the time at which the final counter is removed under
strategy~$S$, and observing that this cannot be less than~$n$, it is
convenient to consider the random variable~$X_S=t_S-n$, taking
non-negative integer values, here denoted the number of ``excess [die]
throws to removal''.

This game is natural and easily specified, and may be played with
nothing more than some plastic counters, a piece of paper, and a pair
of dice.  Applications might include inventory control (for example,
allocating slow-selling goods amongst shops with different sales
rates) and perhaps allocation of multiple indistinguishable tasks
among workers of different efficiencies.

The game is equivalent to a weighted coupon collector's
problem~\cite{boneh1997}, with removing a counter corresponding to
collecting a coupon.  However, in this paper, the allocation of
counters amongst boxes is adjustable by the player at will subject to
the overall sum (quota) being fixed.  The game reduces to that
considered by Myers~\cite{myers2006} if each player allocates the same
number of counters to each box.

\section{The two-box game}

Much of the flavour of the general game is visible in the simplest
non-trivial case: that of two boxes with probabilities~$(p,1-p)$
respectively.  Then the probability mass function
of~$X_{\left<a,b\right>}$, the number of excess throws to removal for
the strategy of placing~$a$ counters in box~1 and~$b$ counters in
box~2, is

\begin{multline}\label{PMF_twoboxes}
  \Prob(\rv{a}{b}=r)=\\ p^a(1-p)^b\left[ {a+b+r-1 \choose a-1,
      b+r}(1-p)^r+
    {a+b+r-1 \choose a+r, b-1}p^r \right],\qquad r=0,1,2,\ldots.
\end{multline}

where~${x+y\choose x,y} = \frac{\left(x+y\right)!}{x!y!}$ is the
choose function.  The two terms correspond to the final counter being
removed from box~1 or box~2 respectively.  Expectations may be
calculated by observing that

\begin{align}
\sum_{r=0}^\infty r {a+b+r-1\choose a+r,b-1}p^r
&= p\sum_{s=0}^\infty\left(s+1\right)p^s{a+b+s\choose a+s+1,b-1}\nonumber\\
&= p\sum_{s=0}^\infty{p^s}\frac{\left(s+1\right)!}{s!}\frac{\left(a+b+s\right)!}{\left(a+s+1\right)!\left(b-1\right)!}\nonumber\\
&= p{a+b\choose a+1, b-1}\sum_{s=0}^\infty
\frac{\Gamma\left(a+b+1+s\right)}{\Gamma\left(a+b+1\right)}
\frac{\Gamma\left(2+s\right)}{\Gamma\left(2\right)}
\frac{\Gamma\left(a+2\right)}{\Gamma\left(a+2+s\right)}
\frac{p^s}{s!}
\nonumber\nonumber\\
&= p{a+b\choose a+1, b-1}\hyper\left(a+b+1,2;a+2;p\right).
\end{align}

\noindent Thus the PMF in Equation~\ref{PMF_twoboxes} has expectation

\begin{multline}\label{expectation}
  \mathbb{E}\left(\rv{a}{b}\right)=
  p^a(1-p)^b\left[
    (1-p){a+b\choose a-1,b+1}\,\hyper\left(a+b+1,2;b+2;1-p\right)+\right.\\
    \left.    p {a+b\choose a+1,b-1}\,\hyper\left(a+b+1,2;a+2;p\right)
    \right]
\end{multline}

where~$\hyper(a,b;c;z)=\sum_{n=0}^\infty\frac{(a)_n
  (b)_n}{(c)_n}\frac{z^n}{n!}$ is the hypergeometric
function~\cite{abramowitz1956}.  Here~$(a)_n=a(a+1)\cdots(a+n-1)$ is
the rising factorial function.  No analytical continuation is needed
here because the primary argument does not intersect the function's
branch cut, conventionally defined as the
interval~$\left[1,\infty\right)$.

\begin{figure}[p]
\includegraphics[width=5in]{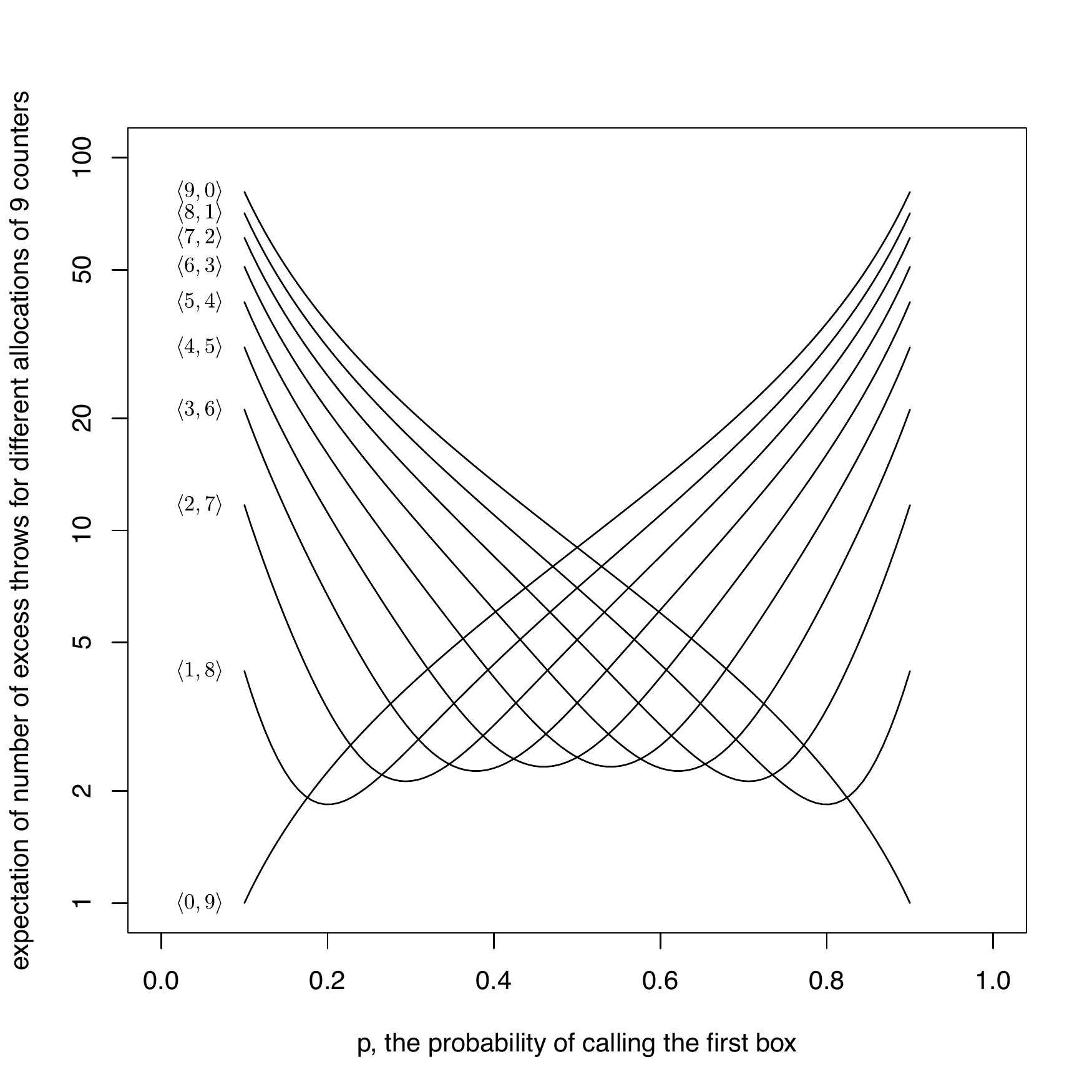} 
\caption{Expected number of excess throws as a function
  of~$p$ \label{expectation_two_categories_highQuality} for different
  allocations of~9 counters in a two-box game; note log scale}
\end{figure}

\begin{figure}[p]
  \includegraphics[width=5in]{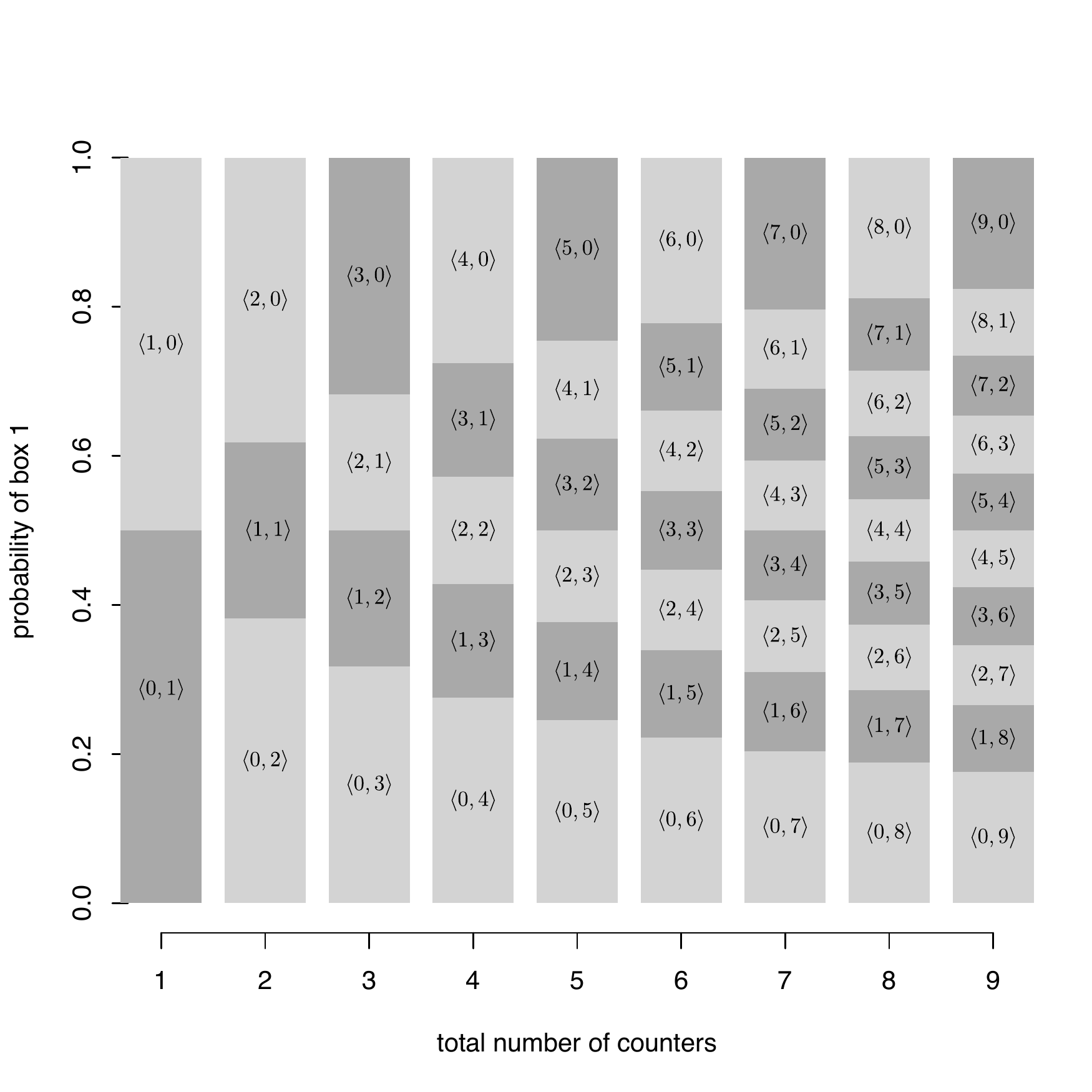} 
  \caption{The two-box game: optimal strategies in expectation
    for~$n=1(1)9$ counters (horizontal
    axis) \label{expectation_cutoff_two_categories} as a function of
    probability of box~1 (vertical axis)}
\end{figure}

The strategy~$\left<a,b\right>$ which minimizes expected time to
removal is
\[\argmin_{a,b|a,b\geqslant 0, a+b=n} \mathbb{E}\rv{a}{b}.\]
Conversely, the cut-off point of~$p$, at which strategy~\strat{a}{b}
yields to strategy~\strat{a-1}{b+1} (assuming all integers
non-negative) is given by solving the equation~$ \mathbb{E}\rv{a}{b}=
\mathbb{E}\rv{a-1}{b+1}$ which may be solved numerically using the
       {\tt hypergeo} package~\cite{hankin2013}.
       Figures~\ref{expectation_two_categories_highQuality}
       and~\ref{expectation_cutoff_two_categories} show some numerical
       results.

\subsection{Game theoretic analysis}

The counter removal process is now considered as a standard two-person
zero sum game with payoff~$+1$ for a victory---that is, removing all
one's counters before the opponent, $-1$ for a loss and~$0$ for a
draw.  The situation is complicated by the fact that play stops when
either player removes her last counter.

There are two natural interpretations: at each time, the players each
randomly choose a box independently from the other players (``separate
[dice] throw''); or alternatively, a random box is chosen for all the
players simultaneously (a ``common throw'').  In the separate throw
case, it is possible to allow the box probabilities to differ between
the players but unless otherwise stated the probabilities will be
identical.

\subsection{Separate throws}

\begin{figure}[p]
  \includegraphics[width=5in]{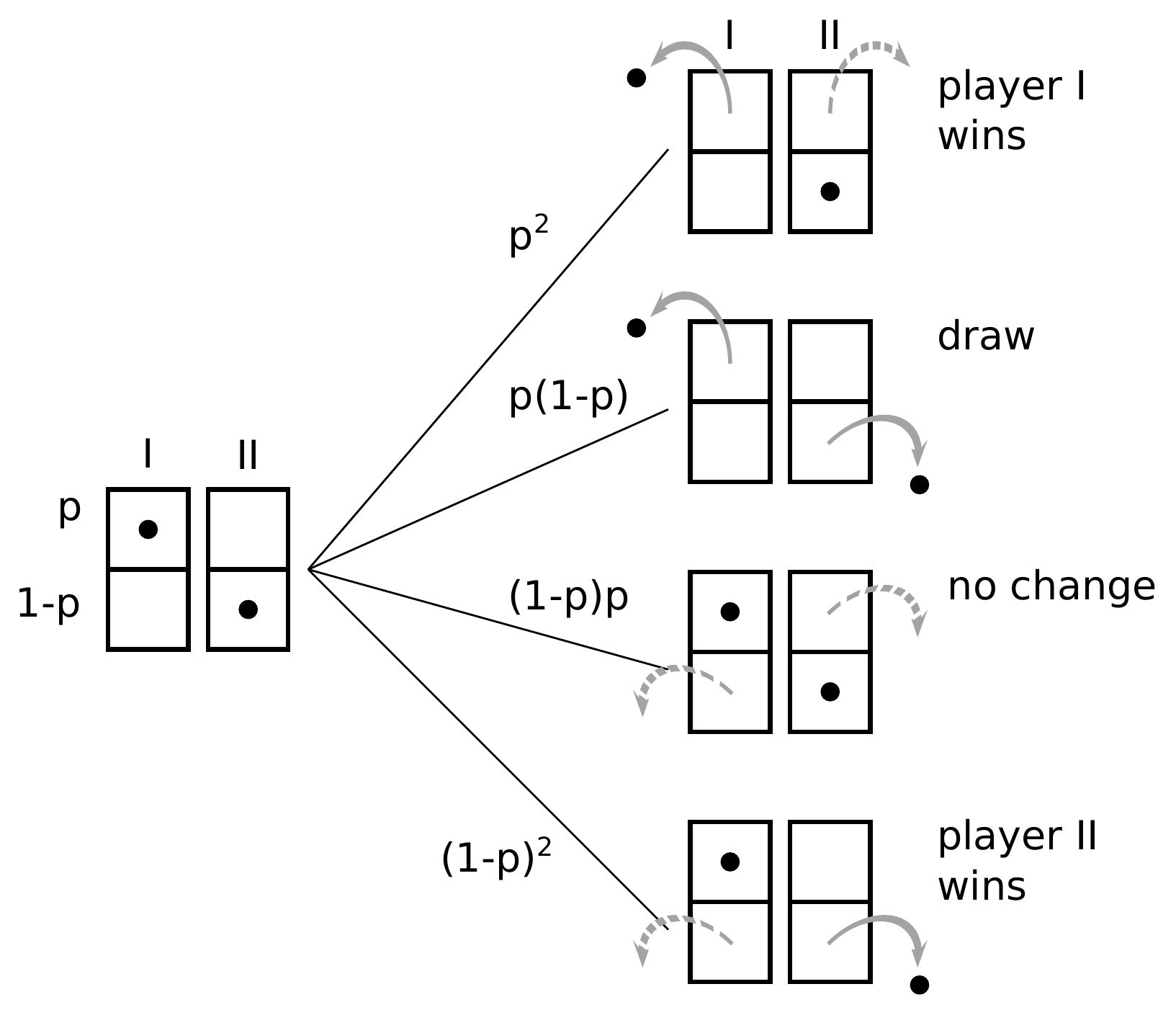}
  \caption{Diagram \label{one_top_one_bottom} showing \strat{1}{0} vs
    \strat{0}{1} game, separate throw game; expected payoff to
    player~I is~$\mathbb{V}$.  Here~$\mathbb{V}=1\cdot p^2 +0\cdot p(1-p) + \mathbb{V}\cdot(1-p)p
    + (-1)\cdot(1-p)^2$, giving~$\mathbb{V}=\frac{2p-1}{1-p(1-p)}$}
\end{figure}

\begin{figure}[p]
  \includegraphics[width=5in]{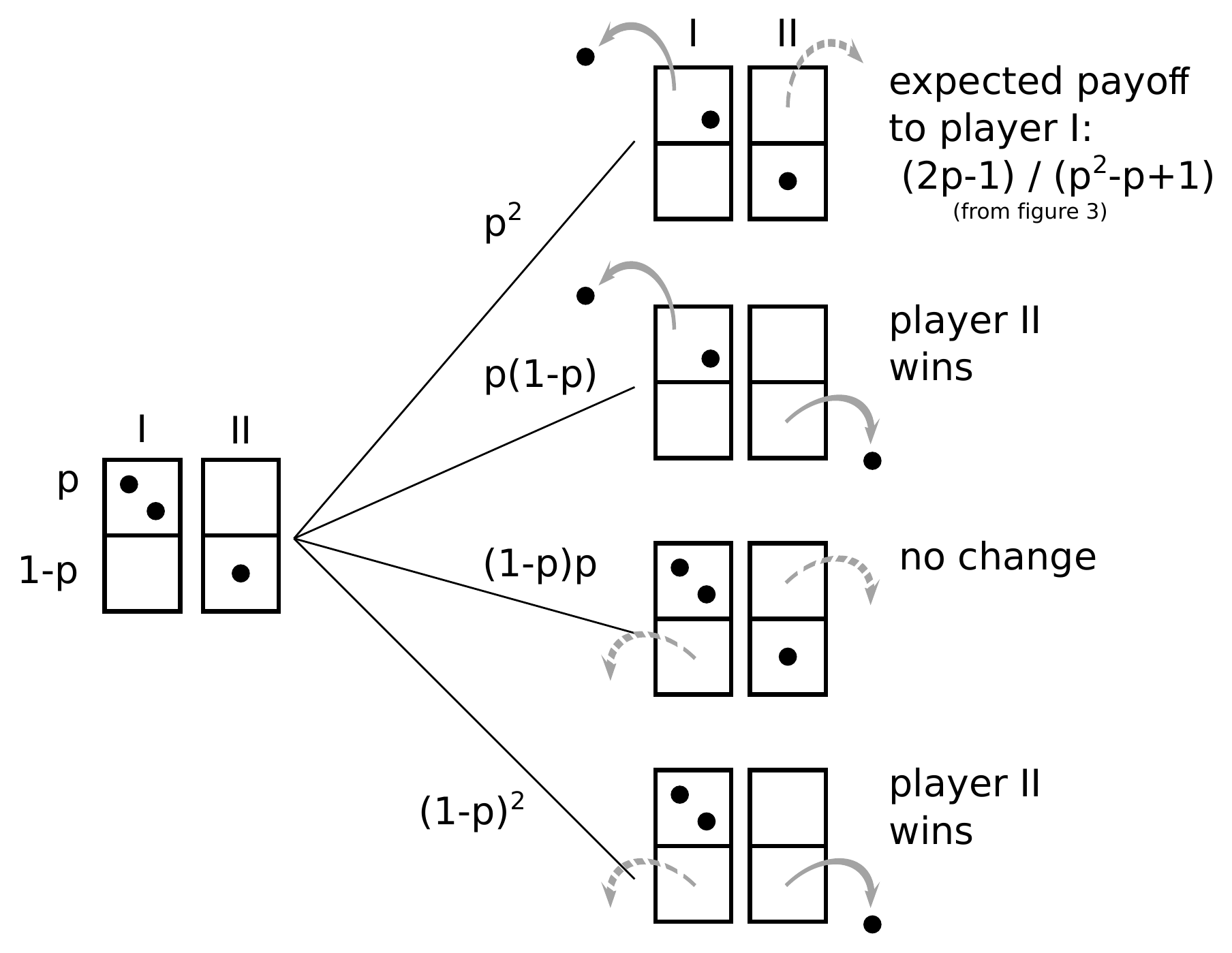}
  \caption{Diagram \label{two_top_one_bottom} showing \strat{2}{0}
    vs \strat{0}{1} game, separate throw game;
    the~$p^2$ term is derived in Figure~\ref{one_top_one_bottom}.  The
    expected payoff to \strat{2}{0}
    is~$\frac{p^2\cdot\frac{2p-1}{1-p(1-p)}+p(1-p)\cdot(-1)+(1-p)^2\cdot(-1)}{1-p(1-p)}=\frac{3p^3-3p^2+2p-1}{(1-p(1-p))^2}$}
\end{figure}

Writing~\ii{a_1}{b_1}{a_2}{b_2} for the game-theoretic payoff to
a~\strat{a_1}{b_1} player versus a~\strat{a_2}{b_2} player,
figures~\ref{one_top_one_bottom} and~\ref{two_top_one_bottom} show
that~$\ii{1}{0}{0}{1}=\frac{2p-1}{p^2-p+1}$
and~$\ii{2}{0}{0}{1}=\frac{3p^3-3p^2+2p-1}{(1-p(1-p))^2}$.  The
simplest non-trivial case with the players possessing an equal number
of counters is~\strat{2}{0} versus~\strat{1}{1} for which the payoff
to the \strat{2}{0} player may be calculated by extending the
methodology of figures~\ref{one_top_one_bottom}
and~\ref{two_top_one_bottom}.  Explicitly, the payoff to
the~$\strat{2}{0}$ player is given as

\begin{align}\label{2011}
\ii{2}{0}{1}{1}   &= p^2\ii{1}{0}{0}{1} \nonumber\\
                   &+ p(1-p)\ii{1}{0}{1}{0}\nonumber\\
                   &+ (1-p)p\ii{2}{0}{0}{1}\nonumber\\
                   &+ (1-p)^2\ii{2}{0}{1}{0}
\end{align}

\begin{figure}[htbp]
\includegraphics[width=5in]{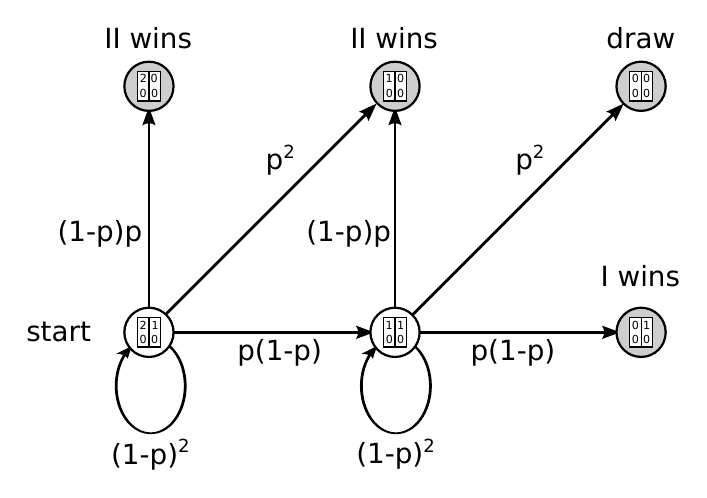}
\caption{Markov chain for~\strat{2}{0} vs~\strat{1}{0}; separate
  calls.  \label{two_top_two_bottom} Absorbing states shown in gray.
  Players~I and~II have box probabilities~$(p,1-p)$.  Using the fact
  that~$\ii{1}{0}{1}{0}=0$, it can be shown that the game value
  is~$1/(p-2)$}
\end{figure}

\begin{figure}[p]
\includegraphics[width=5in]{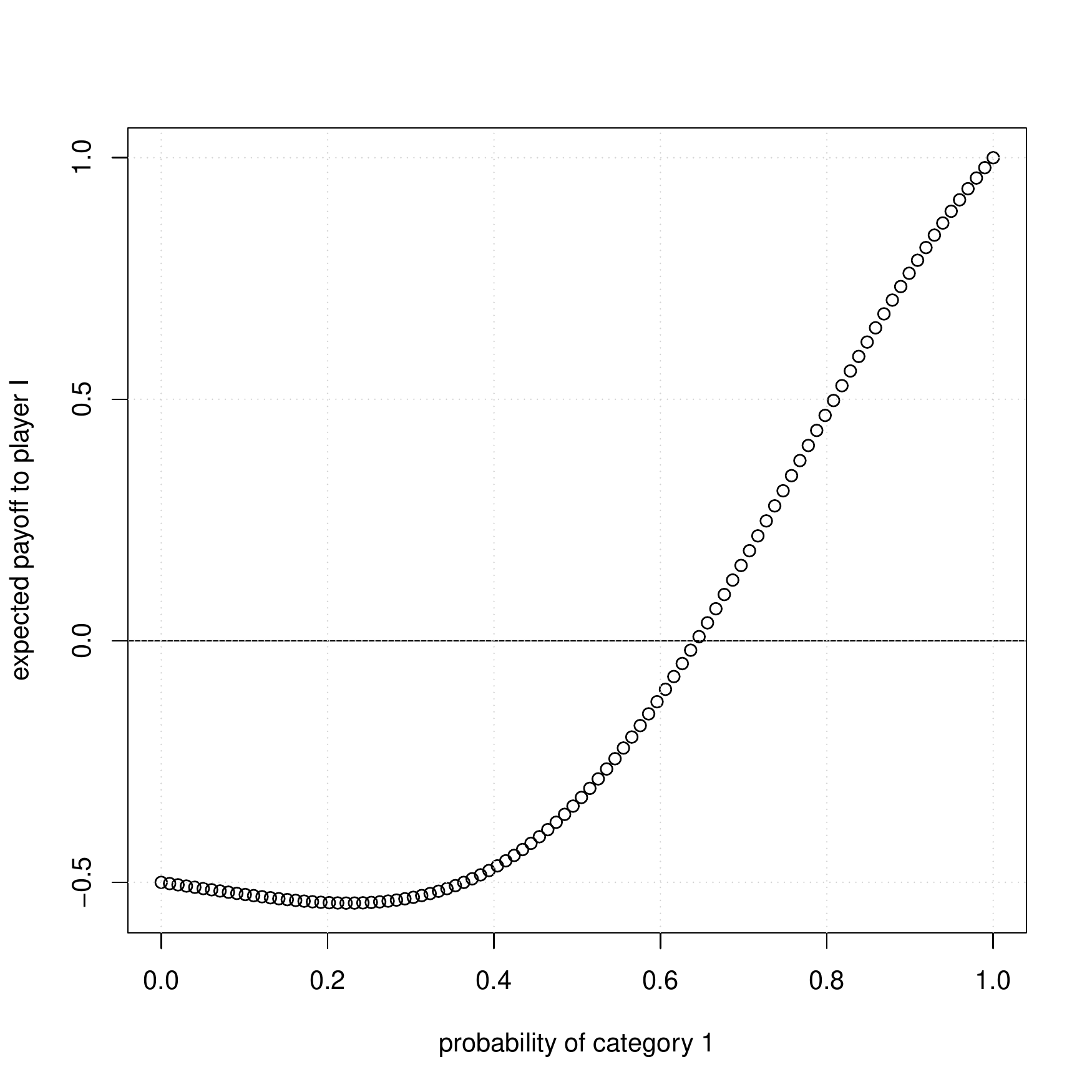}
\caption{Expected payoff to a \strat{2}{0} player against a
  \label{twoboxes_payoff_2011} \strat{1}{1} player as a function
  of (common) box~1 probability; separate calls}
\end{figure}

The first term is given in Figure~\ref{one_top_one_bottom}, the second
is clearly zero, the third given in Figure~\ref{two_top_one_bottom},
and the fourth (in a different format) is given in
Figure~\ref{two_top_two_bottom} as~$1/(p-2)$. Combining these
gives~$\frac{p^5-4p^3+3p^2-2p+1}{p^5-4p^4+7p^3-8p^2+5p-2}$.  Thus
\strat{2}{0} has a positive expected payoff against \strat{1}{1} if
and only if~$p$ exceeds the unique real root of~the numerator, about
0.643.  Figure~\ref{twoboxes_payoff_2011} shows the expected payoff to
\strat{2}{0} as a function of~$p$ and exhibits an unexpected feature:
a nontrivial local minimum near~$p=0.225$ of
about~$\mathbb{V}=-0.543$.  This observation is not useful in the
general context of this paper in which a specified number of counters
must be allocated, but might be relevant if a player is {\em
  constrained} to play~\strat{1}{1} against~\strat{2}{0} (and may
vary~$p$ at will).

\begin{figure}[p]
\includegraphics[width=5in]{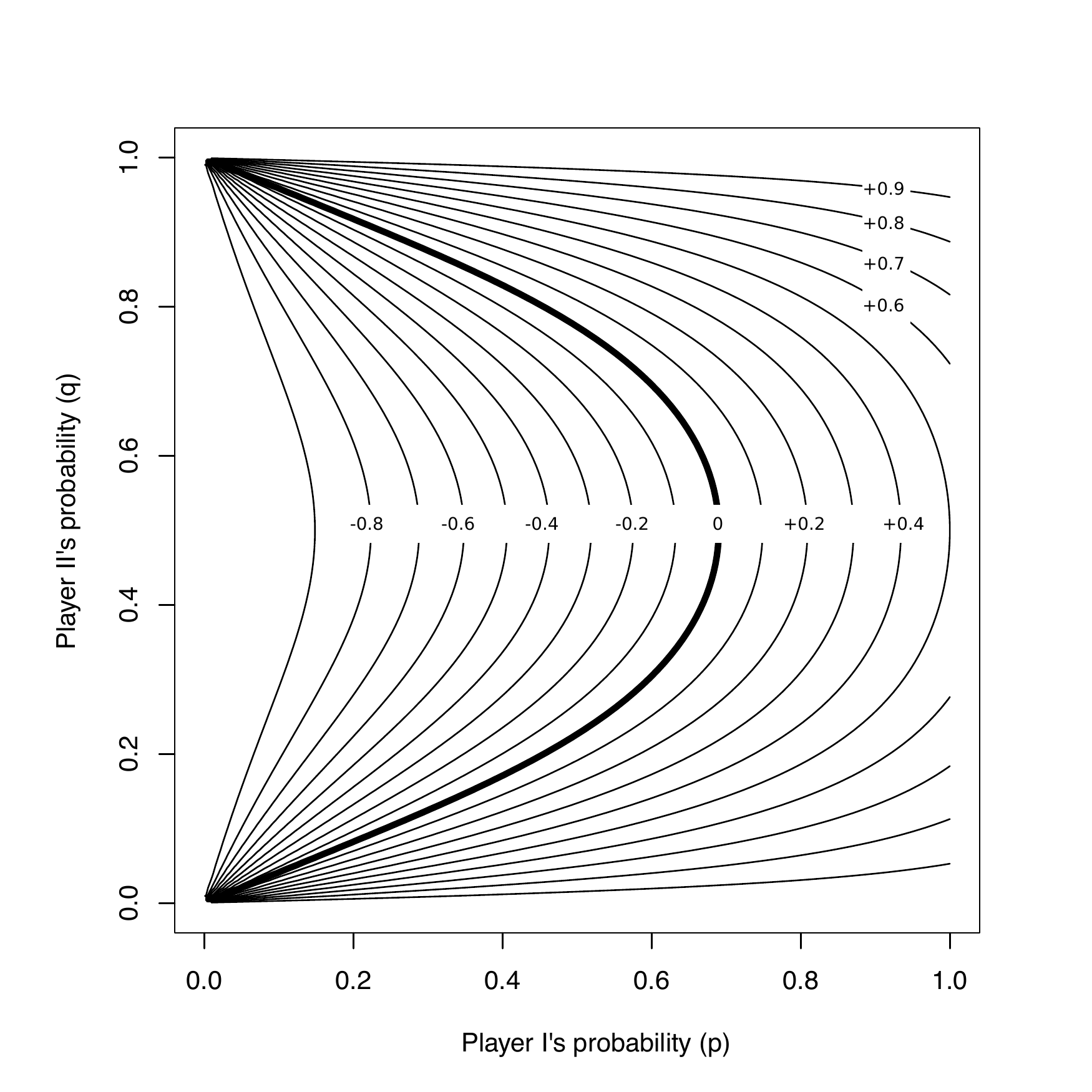} 
\caption{Expected payoff to \strat{2}{0} (box probablities~$p,1-p$)
 \label{contour_pq} against \strat{1}{1} (box probablities~$q,1-q$)}
\end{figure}

It is possible to relax the requirement that the players'
probabilities are identical, although the algebra becomes complicated
quickly.  Figure~\ref{contour_pq} shows the expected payoff to I as a
function of the players' probabilities of choosing the first box.
If viewed as a game on the unit square in the sense of
Heuer~\cite{heuer2001}, the expected payoff has a saddle point
at~$p=1,q=1/2$ at which the payoff to~I is 0.5.

%

\subsection{Two boxes, common throw}

\begin{figure}[p]
\includegraphics[width=5in]{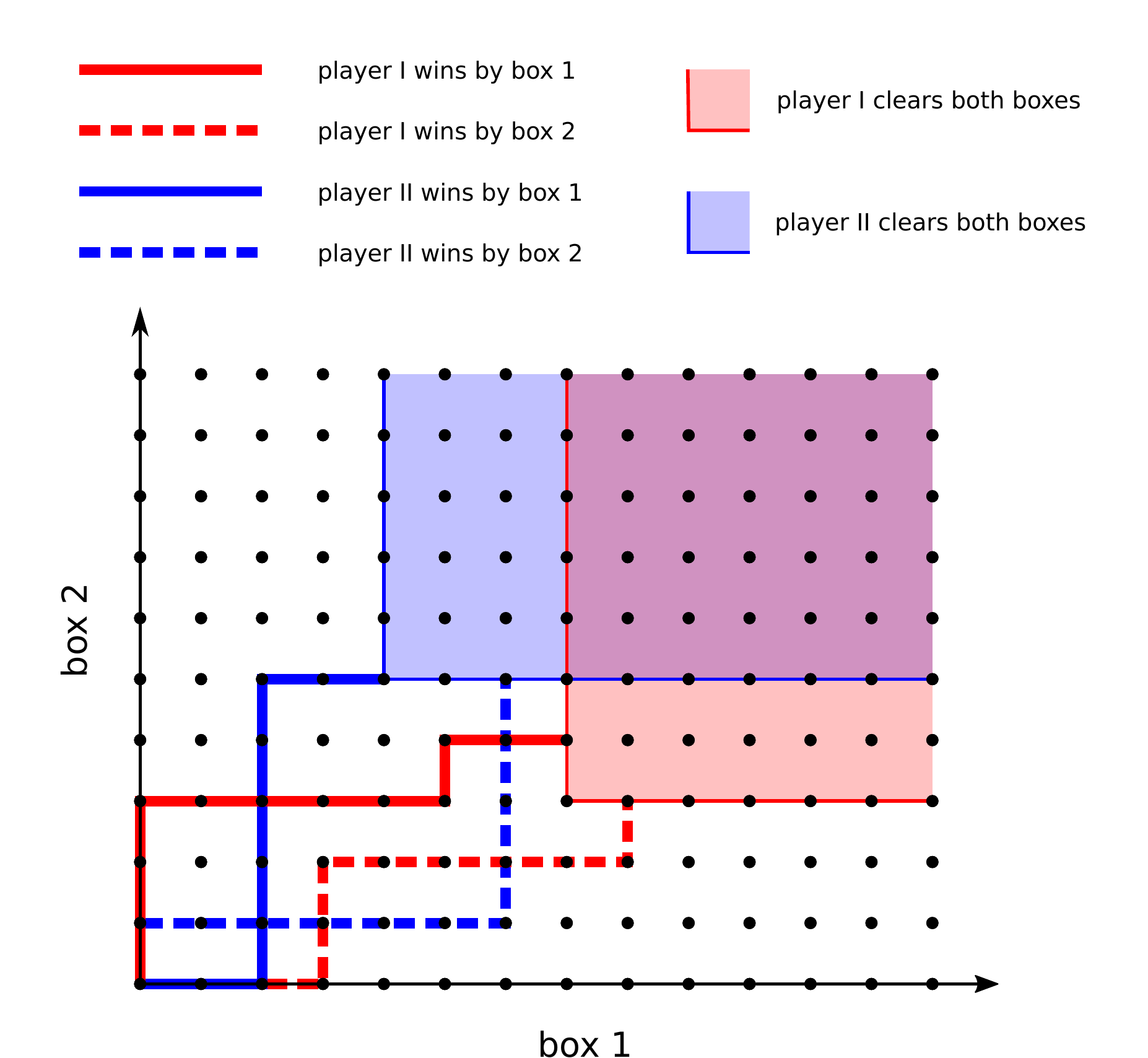} 
\caption{Two box, common thow game with \label{path_twoboxes} player~I
  playing~\strat{8}{3} and player~II playing~\strat{4}{5}.  A game is
  represented by a coloured line in the form of a random walk from the
  origin; axes indicate total number of box 1 and box 2 throws.  Line
  colours show the winner: red for player I winning and blue for
  player II winning, and line type shows the winning box: solid for
  winning by box 1 and dotted for winning by box 2.  Shaded regions
  indicate a player clearing both boxes: red for player I and blue for
  player II}
\end{figure}

The case where the throw is common to both players is qualitatively
different; see Figure~\ref{path_twoboxes}.  In the following, we
suppose that~$a,b,m,n$ are strictly positive integers.

Suppose player~I plays~\strat{a+n}{b} and player~II
plays~\strat{a}{b+m}.  Then player I may win in one of two ways: the
final counter to be removed may be from box~1 (written ``winning by
box 1"), or from box~2.

The successive throws constitute an iid sequence of~1's and~2's.  At
any time, denote the total number of box~1 throws as~$r_1$ and number
of box~2 throws as~$r_2$.  Player~I wins by box~2 if, when~$r_2=b$, we
have $r_1\geqslant a+n$, and this occurs with
probability~$\Prob(r_1\geqslant a+n|r_2=b)$, a negative binomial
distribution.  Similarly, player I wins by box~1 with
probability~$\Prob(b\leqslant r_2< b+m|r_1=a+n)$.

Because the complementary cumulative distribution function of the
negative binomial with parameters~$r,p$ is given by the regularized
incomplete beta function\footnote{This is a standard function defined
  as~$I_x(\alpha,\beta)=\frac{\Gamma(\alpha+\beta)}{\Gamma(\alpha)\Gamma(\beta)}\int_{p=0}^x
  p^{\alpha-1}(1-p)^{\beta-1}\,dp$.  Here it is evaluated using the
  GSL library~\cite{hankin2006a} or the hypergeo
  library~\cite{hankin2013}.}~$I_p(k+1,r)$, the probabilities of
player I beating player II by box 1
is~$1-I_p\left(a+n,b+m\right)-\left(1-I_{1-p}\left(b,a+n\right)\right)$
and by box~2 is~$I_{p}\left(a+n,b\right)$.

Thus the total probability of player~I winning
is~$I_{p}\left(a+n,b+m\right)$ and,\\ observing that a draw is
impossible, the expected payoff to player~I is
thus\\ \mbox{$2I_{1-p}\left(b+m,a+n\right)-1$}.  Note that these
probabilities are functions only of the larger number of counters in
each of the two boxes.

Similar considerations apply if player~I plays~\strat{a}{b} and~II
plays~\strat{a+m}{b}; now player~II cannot win but draws with
probability~$I_{1-p}\left(b+m,a\right)$ (both players removing their
final counter from box~2 simultaneously).

For completeness, we observe that~\strat{a}{b} certainly
beats~\strat{a+m}{b+n}, and further that~\strat{a}{b} certainly draws
against~\strat{a}{b}, even if one or both of~$a$ and~$b$ is zero.

\begin{figure}[p]
\includegraphics[width=5in]{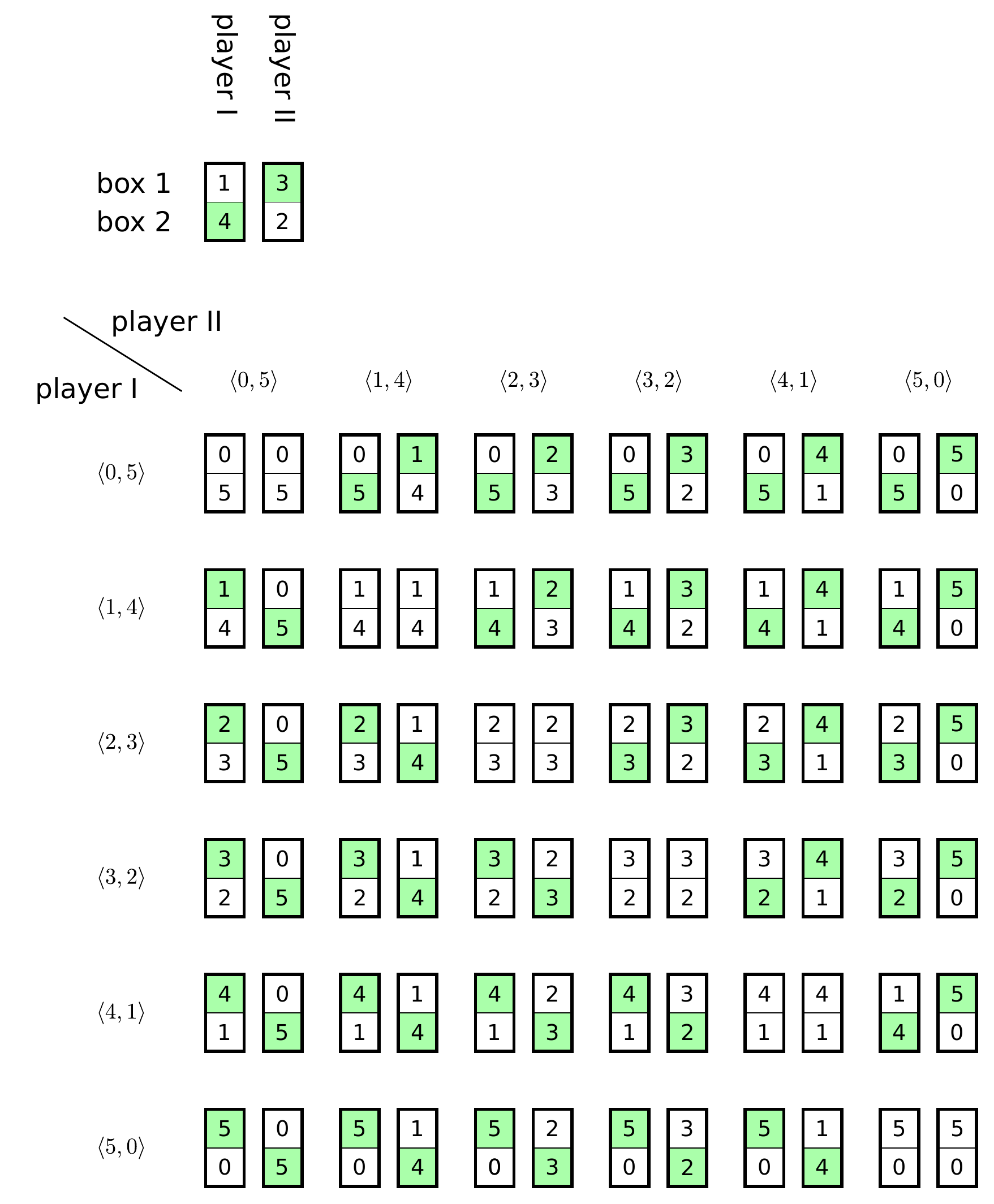} 
\caption{All~$6\times 6=36$ possible games played between two players
  with five counters each.  For each game, the two rectangles
  represent the\label{twobox_fits_on_page} allocation of counters by
  player~I (left) and player~II (right); green indicates that that
  box's counters are ``active'' in the sense that the number appears
  in the expression for the payoff to player~I,
  viz~$2I_{1-p}\left(b+m,a+n\right)-1$}
\end{figure}

The minimax strategy for the common-call two-category game is pure,
except possibly at boundaries.  This follows from the facts that: the
payoff matrix is antisymmetric (see Figure~\ref{twobox_fits_on_page}
and Table~\ref{payoff_matrix_5counters_p=0.4}); each row is decreasing
with column number on the left of the diagonal, and increasing to the
right (because $I_p(a-1,b)>I_p(a,b)$ for~$a\geqslant 1$); and, for
any~$p\in\left[0,1\right]$ one has~$I_p(a+1,b)>I_p(a,b+1)$.

Critical points for~$p$ are found by solving~$I_p(r,n+1-r)=1/2$
for~$r=0,1,\ldots,n$.  For example, the range of~$p$ for which
\strat{1}{4} is the dominant strategy~$\left(p_L,p_U\right)$ is given
by~$I_{p_L}(4,2)=I_{p_U}(5,1)=1/2$; numerically this is about~$(0.686,
0.871)$.

\section{Three boxes}

The three box case is now considered, partly motivated by the fact
that the space of probabilities is easily visualized on a ternary
plot.

\begin{figure}[p]
  \includegraphics[width=5in]{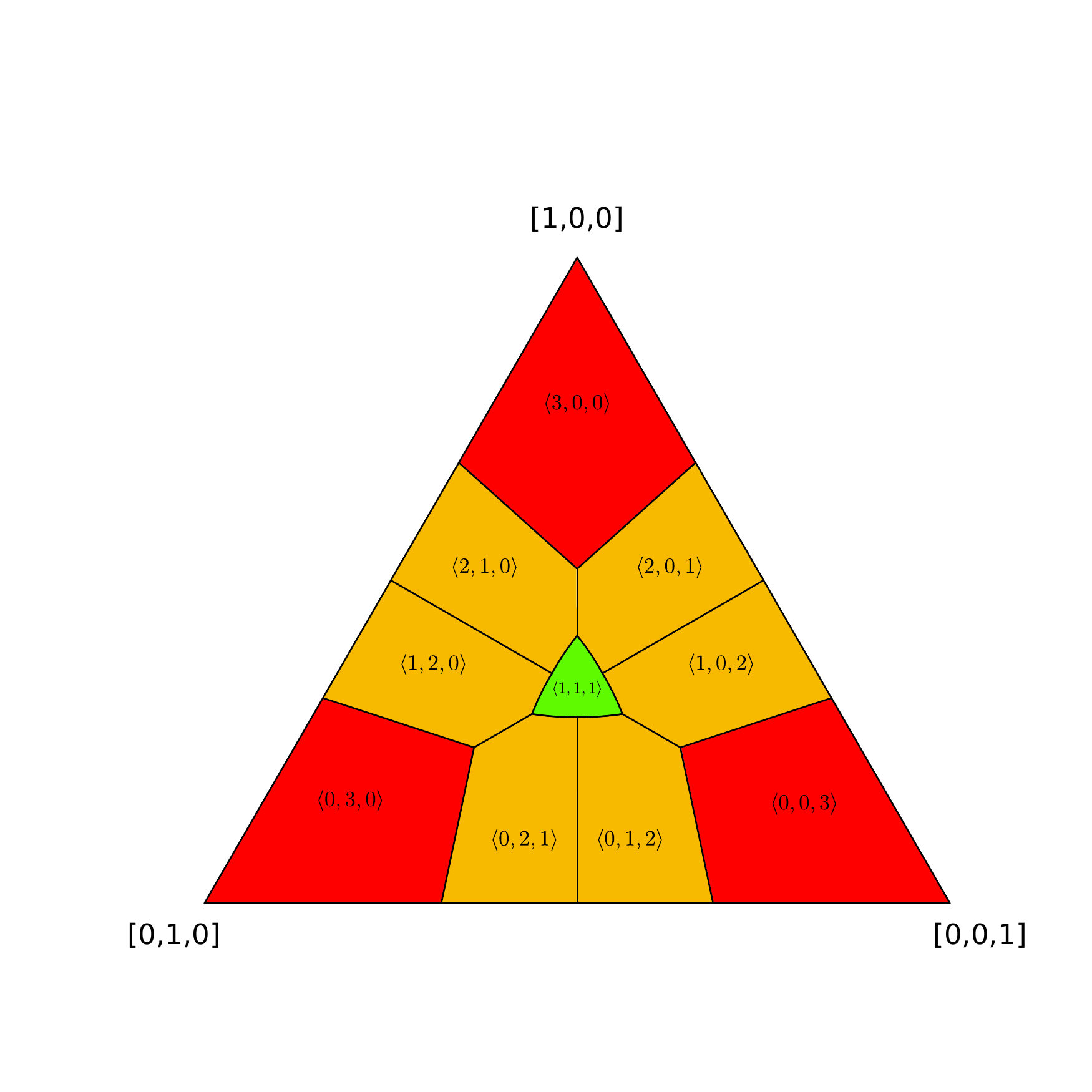}
  \caption{Triangular
    plot \label{3boxes_3counters_expectation_analytic_post_inkscape}
    showing optimal strategies in expectation as a function of the
    three box probabilities~$p_1,p_2,p_3$, where~$p_1+p_2+p_3\leq
    1$. Different coloured regions show the extent of the different
    strategies}
\end{figure}

The simplest nontrivial case for three boxes is with three counters,
in which case there are are 10 strategies, from~$\left<3,0,0\right>$
through~$\left<1,1,1\right>$ to~$\left<0,0,3\right>$. The optimal
strategy for any box probabilities is shown in
figure~\ref{3boxes_3counters_expectation_analytic_post_inkscape}.

\begin{figure}[p]
  \includegraphics[width=5in]{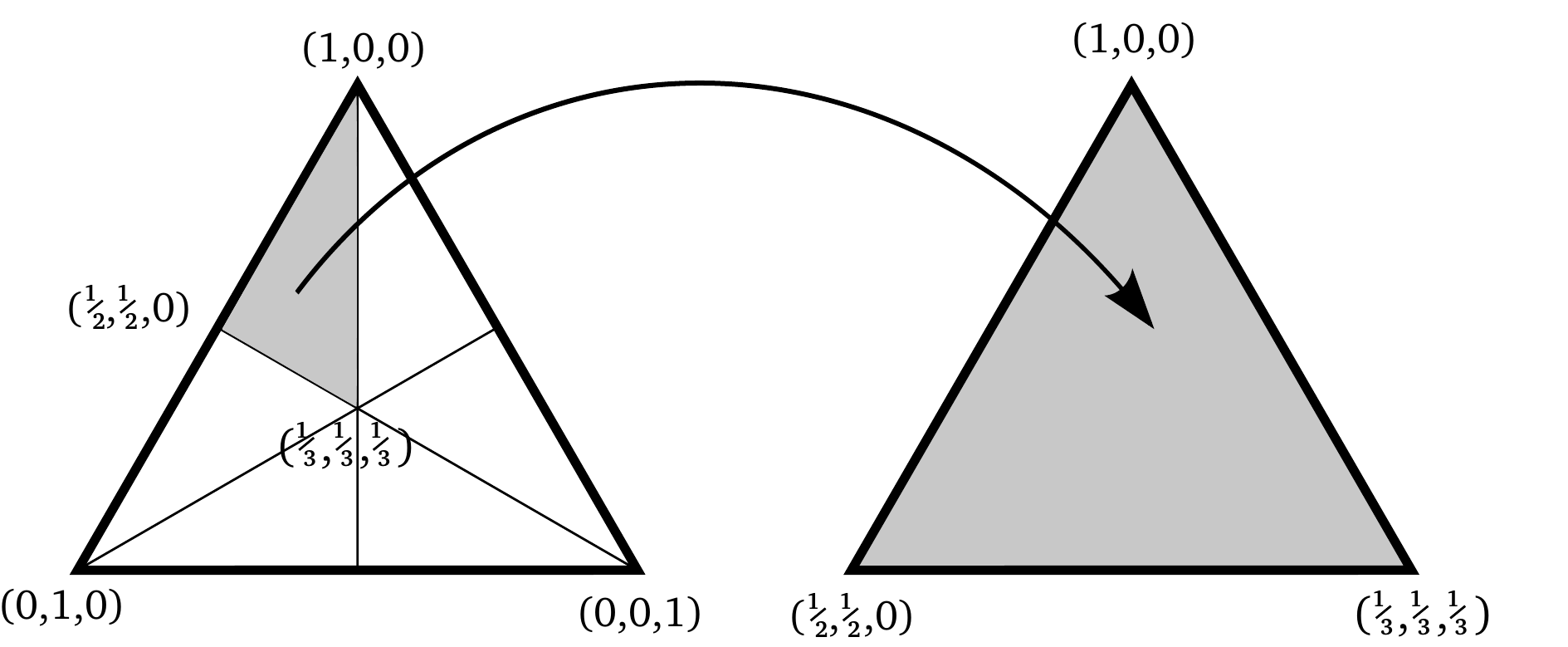}
  \caption{Diagram showing mapping used for the 
    \label{ternary_diagram_revised} ternary diagrams}
\end{figure}

It is possible to streamline graphical results by taking advantage of
the fact that we may assumet, without loss of generality, that
$p_1\geqslant p_2\geqslant p_3$
(Figure~\ref{ternary_diagram_revised}).  In what follows, we will
consider only strategies for which~$n_1\geqslant\ n_2\geqslant
n_3$---restricted partitions in the language of
Hankin~\cite{hankin2006}---as strategies not satisfying this
constraint may be improved in expectation\footnote{To see this,
  consider the general case of~$p_1\geqslant\cdots\geqslant p_k$ and a
  strategy that violates the constraint, that is, there exist~$i,j$
  with~$1\leqslant i<j\leqslant k$ and~$n_i<n_j$.  Then split the
  boxes into two sets:~$S=\left\{i,j\right\}$ and its
  complement~$\overline{S}$.  Then removal of a counter from either
  set may be embedded into two independent Poisson processes,
  following~\cite{boneh1997}.  Then consider transferring a counter
  from box~$i$ to box~$j$.  The~$\overline{S}$ process is unaffected,
  but the~$S$ process has lower expected time to completion.  Proof:
  consider the common throw game with player I playing~\strat{a+1}{b}
  and player II \strat{a}{b+1}; here $a>b$ and~$p>\frac{1}{2}$.  We
  shall show that player I's expected throws is strictly less than
  player II's.  First, define I's margin of victory~$M$ as the number
  of throws needed for II to clear her counters after I clears her
  counters ($M<0$ means she lost; note that~$M\neq 0$).
  Then~$\mathbb{E}(M)$ may be calculated
  as~$\mathbb{E}\left(M\left|M>0\right.\right)\cdot\Prob(M>0)+
  \mathbb{E}\left(M\left|M<0\right.\right)\cdot\Prob(M<0)$.

  From Figure~\ref{path_twoboxes}, we see that the two conditional
  distributions are geometric and
  so~$\mathbb{E}(M)=\frac{I_p\left(a+1,b+1\right)}{1-p} -
  \frac{1-I_{p}\left(a+1,b+1\right)}{p}$.  This is strictly positive
  because~$1-p<\frac{1}{2}<I_{1/2}\left(a+1,b+1\right)<
  I_{p}\left(a+1,b+1\right)$.  Thus a player may increase her game
  value by transferring a counter from box~$i$ to box~$j$.}  and game
value\footnote{Proof: for the common throw case, we observe
  that~$p>\frac{1}{2}$ and~$a>b$
  imply~$I_{p}\left(a+1,b+1\right)>\frac{1}{2}$.  For the separate
  throw case we observe that the existence of~$a,m,n>0$ with
  $p>\frac{1}{2}$ and~$\payoff{\strat{a+n}{a}}{\strat{a}{a+m}}<0$
  implies that $\exists a'<a$ and~$m',n'>0$
  with~$\payoff{\strat{a'+n'}{a}}{\strat{a'}{a'+m'}}<0$.  The result
  follows from the fact
  thatthat~$\payoff{\strat{0}{0}}{\strat{0}{0}}=0$.}.

\subsection{Expectation}

The probability mass function for~$X_{\left<n_1,n_2,n_3\right>}$ (that
is, the time to removal of the final counter, minus~$n=n_1+n_2+n_3$) is
given by

\begin{align}
  \Prob(X_{\left<n_1,n_2,n_3\right>}=r) &=
       \sum_{r_2+r_3=r} {n+r_2+r_3-1\choose n_1-1  ,n_2+r_2,n_3+r_3}p_1^{n_1    } p_2^{n_2+r_2} p_3^{n_3+r_3}\nonumber\\
    &+ \sum_{r_1+r_3=r} {n+r_1+r_3-1\choose n_1+r_1,n_2-1  ,n_3+r_3}p_1^{n_1+r_1} p_2^{n_2    } p_3^{n_3+r_3}\nonumber\\
    &+ \sum_{r_1+r_2=r} {n+r_1+r_2-1\choose n_1+n_1,n_2+r_2,n_3-1  }p_1^{n_1+r_1} p_2^{n_2+r_2} p_3^{n_3    }.\label{threesquares}
  \end{align}

Here, the three terms correspond to the final counter being removed
from boxes~1,2,3 respectively; and~$r_i$ corresponds to the number of
``wasted'' box~$i$ throws (that is, the number of times box~$i$ is
called when box~$i$ is empty).  This equation is analytically
challenging; the first term is algebraically equal to

\begin{multline}\label{ithink1}
  p_1^{n_1} p_2^{n_3}  p_3^{n_3}\left[
  {p_3}^r{n+r-1\choose n_1-1, n_2, n_3+r}\hyper\left(1,-n_3-r;n_2+1;-\frac{p_2}{p_3}\right)\right. -\\
  \left.
  {p_2}^{r}\frac{p_2}{p_3}{n+r-1\choose n_1-1,n_2+r+1,n_3-1}\hyper\left(1,1-n_3;n_2+r+2;-\frac{p_2}{p_3}\right)
  \right].
\end{multline}

However, because the hypergeometric functions each have an upper
argument of a strictly negative integer, the expressions are a
polynomial in~$p_2/p_3$.  However, equation~\ref{ithink1} is the
preferred form for numerical evaluation~\cite{hankin2013}.  The other
terms are analogous, unless one of the~$n_i=0$, in which case that
term is omitted.  These results may be used to show that the expected
times to removal for the three consistent strategies are:

\begin{align}
\mathbb{E}X_{\left<3,0,0\right>} &= \frac{3(1-p_1)}{p_1}\\
\mathbb{E}X_{\left<2,1,0\right>} &= -3+\frac{1}{p_2} + \frac{p_2(3p_1+2p_2)}{p_1(p_1+p_2)^2}\\
\mathbb{E}X_{\left<1,1,1\right>} &=
-2+
\frac{1}{p_1}+
\frac{1}{p_2}+
\frac{1}{p_3}
-\frac{1}{1-p_1}
-\frac{1}{1-p_2}
-\frac{1}{1-p_3}
\end{align}

(the first result is more easily determined as the expectation of an
inverse binomial distribution with parameters 2 and~$p_1$, and the
third is a special case of the Coupon Collector's problem
with~$n=3$~\cite{berenbrink2009}).

\begin{figure}[p]
  \includegraphics[width=5in]{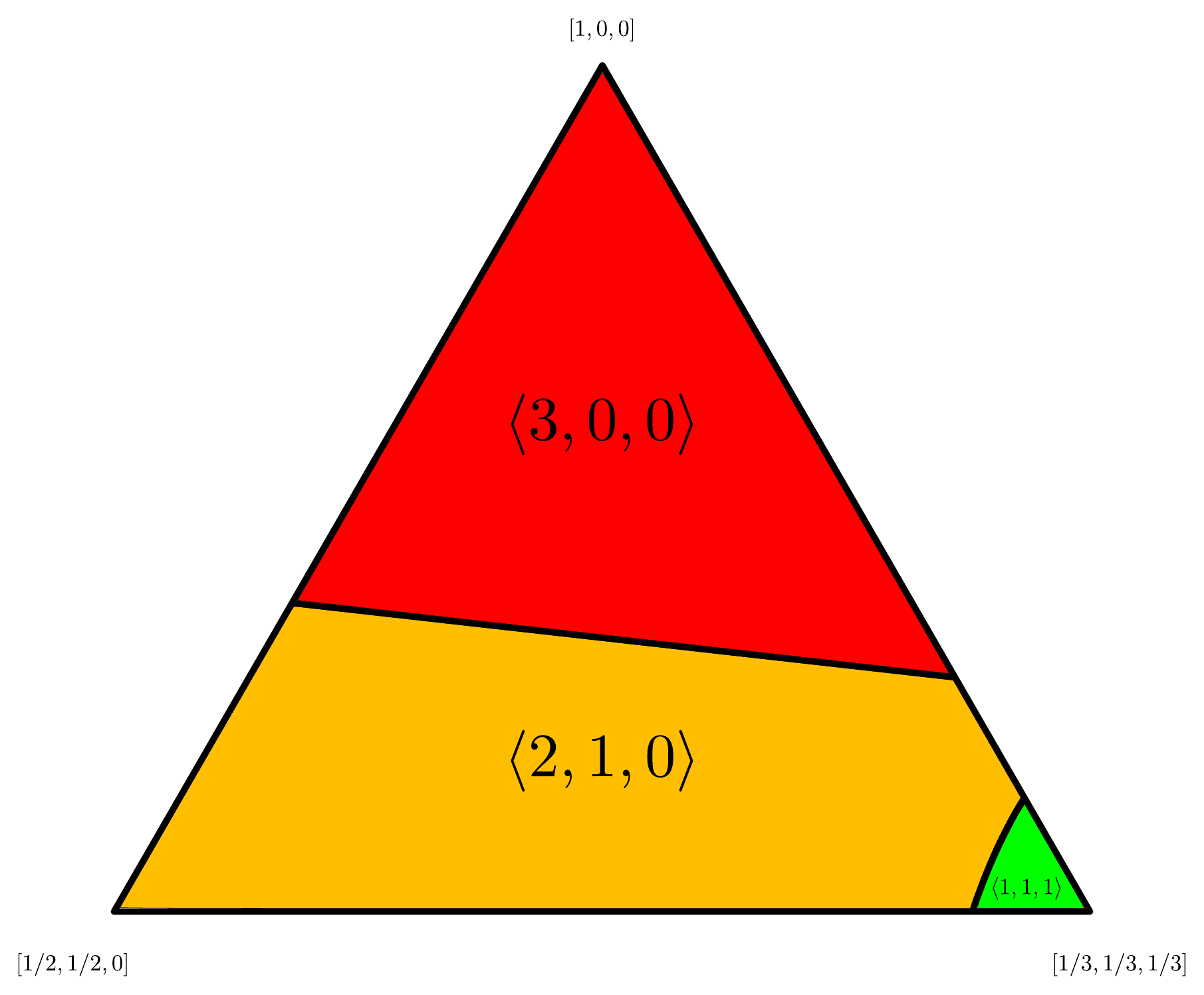}
  \caption{Triangular plot \label{best_expectation_3counters_3boxes}
    showing optimal strategy in expectation for allocating 3 counters
    among 3 boxes, as a function of~$p_1,p_2,p_3;p_1\geqslant
    p_2\geqslant p_3$.  The three points of the triangle correspond
    to~$\left(1,0,0\right)$, $\left(\frac{1}{2},\frac{1}{2},0\right)$,
    and~$\left(\frac{1}{3},\frac{1}{3},\frac{1}{3}\right)$.  The line
    common to the~$\left<3,0,0\right>$ and~$\left<2,1,0\right>$ regions is
    straight, but the line common to~$\left<2,1,0\right>$
    and~$\left<1,1,1\right>$ is not}
\end{figure}

To find the boundaries of the regions in which the three strategies
are optimal in expectation,
solve~$\mathbb{E}X_{\left<3,0,0\right>}=\mathbb{E}X_{\left<2,1,0\right>}$
and~$\mathbb{E}X_{\left<2,1,0\right>}=\mathbb{E}X_{\left<1,1,1\right>}$.
These yield relationships showing that while the boundary of the
$\left<3,0,0\right>$ region is linear, the boundary
of~$\left<1,1,1\right>$ is not; thus not all the regions in
Figure~\ref{best_expectation_3counters_3boxes} are polygonal, despite
appearances.

\subsubsection{More than three counters}

For the three-box case with an arbitrary number of counters, one
natural strategy would be to allocate one's quota in proportion
to~$p_1,p_2,p_3$ (or at least as close as this as possible).  However,
this is not necessarily optimal in expectation.  Consider the game
with probabilities~$(0.7,0.2,0.1)$ and ten counters; one might think
strategy~$\left<7,2,1\right>$ to be superior to, say,
$\left<9,1,0\right>$.  However, numerical methods show that the former
strategy has expected time to removal of~$6.09$ (to 2dp), compared
with the latter which has~$3.37$.

\begin{figure}[p]
  \includegraphics[width=5in]{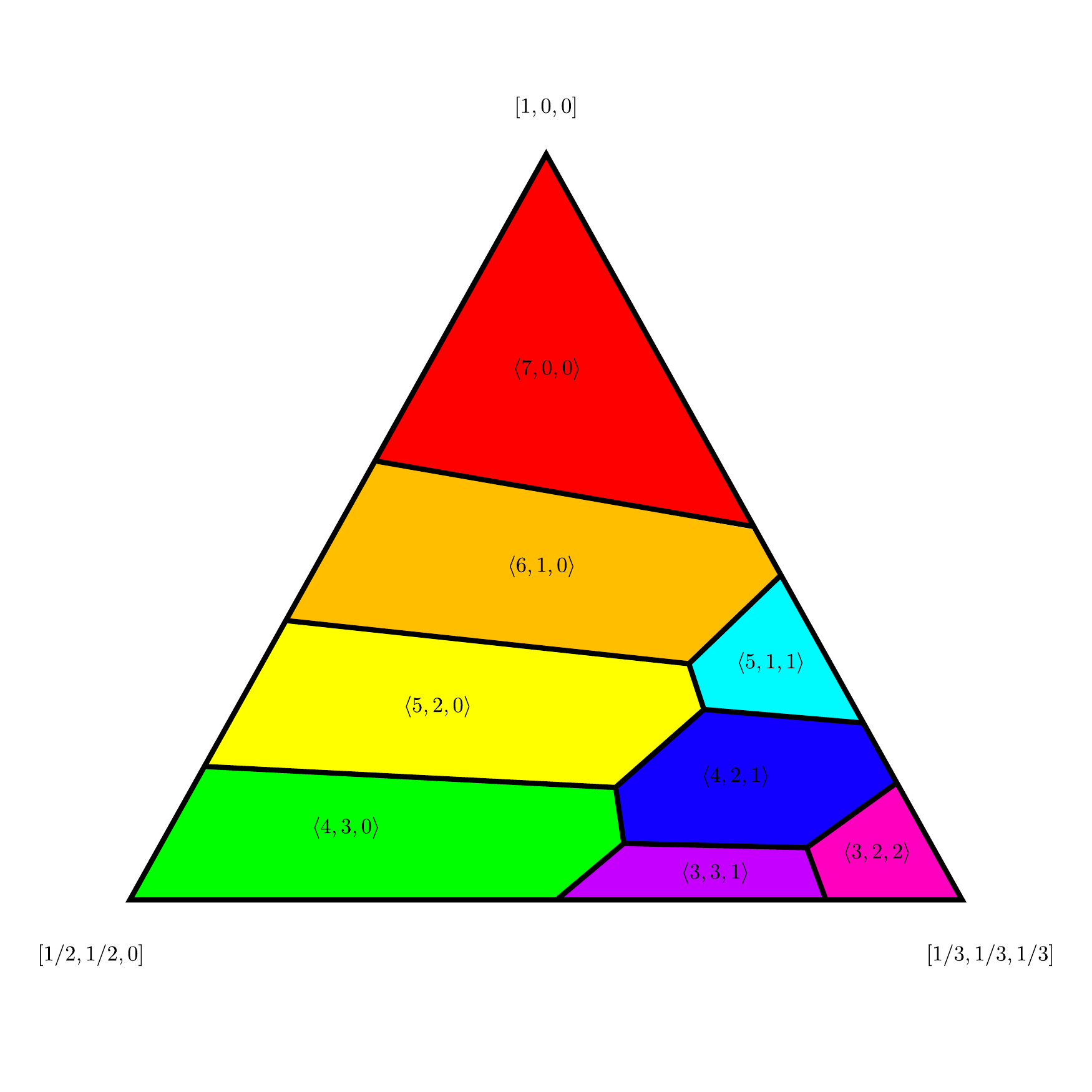}  
  \caption{Triangular plot \label{whichisthebest} showing optimal
    strategy in expectation for allocating 7 counters among 3 boxes,
    as a function of~$p_1,p_2,p_3, p_1\geqslant p_2\geqslant p_3$}
\end{figure}

Figure~\ref{whichisthebest} shows a triplot of the best strategy in
expectation for distributing seven counters among three boxes over the
region~$p_1>p_2>p_3$.

\subsection{Game theoretic analysis}

For both the common throw game and the separate throw game, natural
extensions of the methods of Figures~\ref{one_top_one_bottom}
and~\ref{two_top_one_bottom} may be applied.  Here the simplest non
trivial case, that of three counters, is considered.  For the common
throw game, Table~\ref{payoff_matrix_3_boxes_3counters_common} shows
the payoff matrix as a function of the probabilities; but the
equivalent calculation for the separate case is algebraically more
involved.  To illustrate of this, if
Table~\ref{numeric_payoff_matrix_3counters_3_boxes_p=0.6_0.3_0.1} is
written in exact, rational form with fractions in their lowest terms,
one of the entries has denominator exceeding~$5\times 10^{27}$.  In
the absence of intelligible algebraic formulas,
Table~\ref{numeric_payoff_matrix_3counters_3_boxes_p=0.6_0.3_0.1}
shows the payoff matrix for the three counters, three box case,
written to three decimal places, for both the common throw and
separate throw games.  The tables show perhaps unexpected differences;
for example, the {\em sign} of the payoff to~$\left<3,0,0\right>$ when
playing against~$\left<2,1,0\right>$ changes; also note that the common
throw game has a saddle point at~$\left<2,1,0\right>$, the separate
throw game at~$\left<3,0,0\right>$.

\begin{figure}[p]
\includegraphics[width=5in]{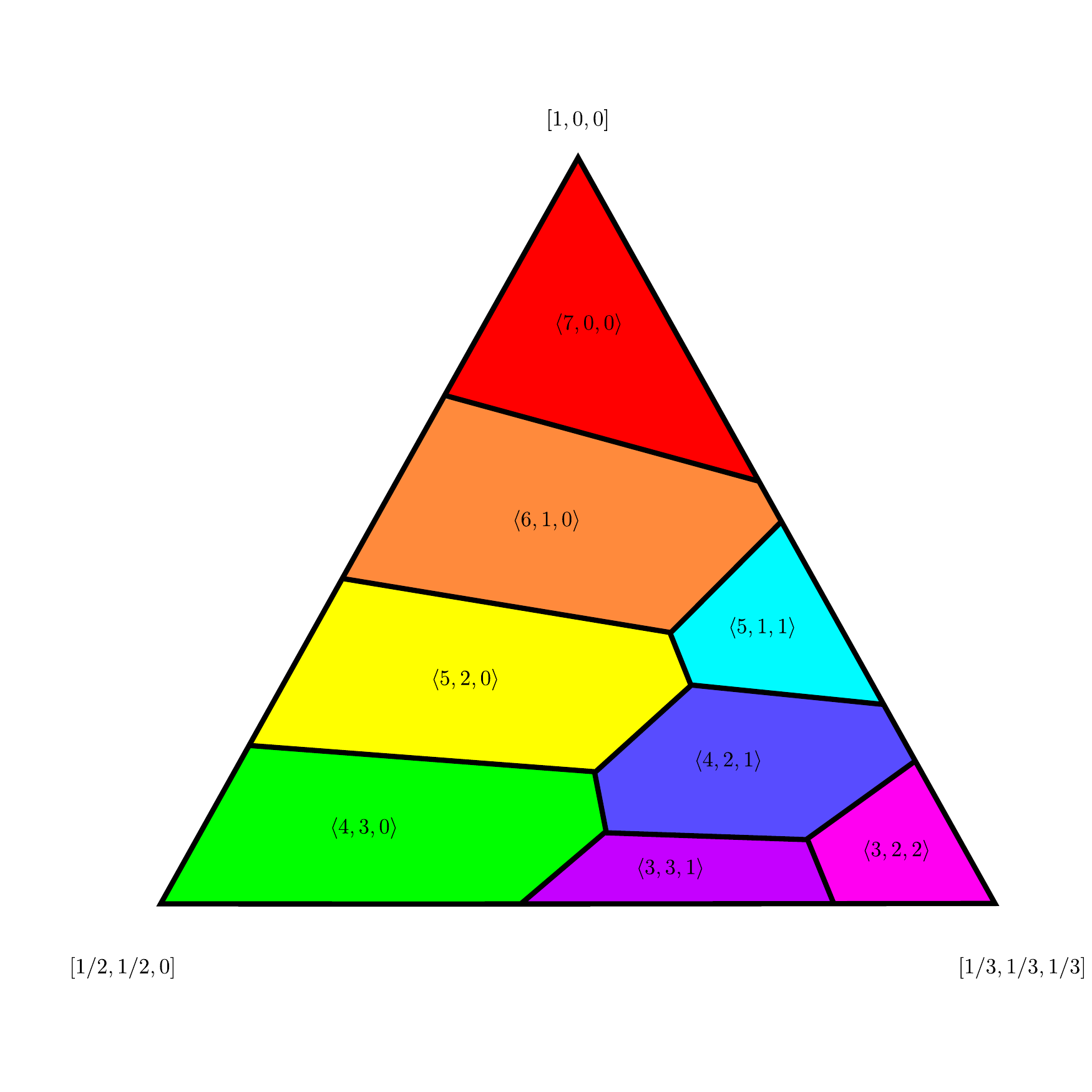}  
\caption{Triangular plot \label{winner_7sep} showing regions which are
  game-theoretic optimal for the 7 counter, 3 box case (separate
  throws)}
\end{figure}

It is possible to conduct numerical experiments for larger numbers of
counters.  Figure~\ref{winner_7sep} shows the appropriate strategy for
the case with three boxes and seven counters.

\begin{figure}[p]
  \includegraphics[width=5in]{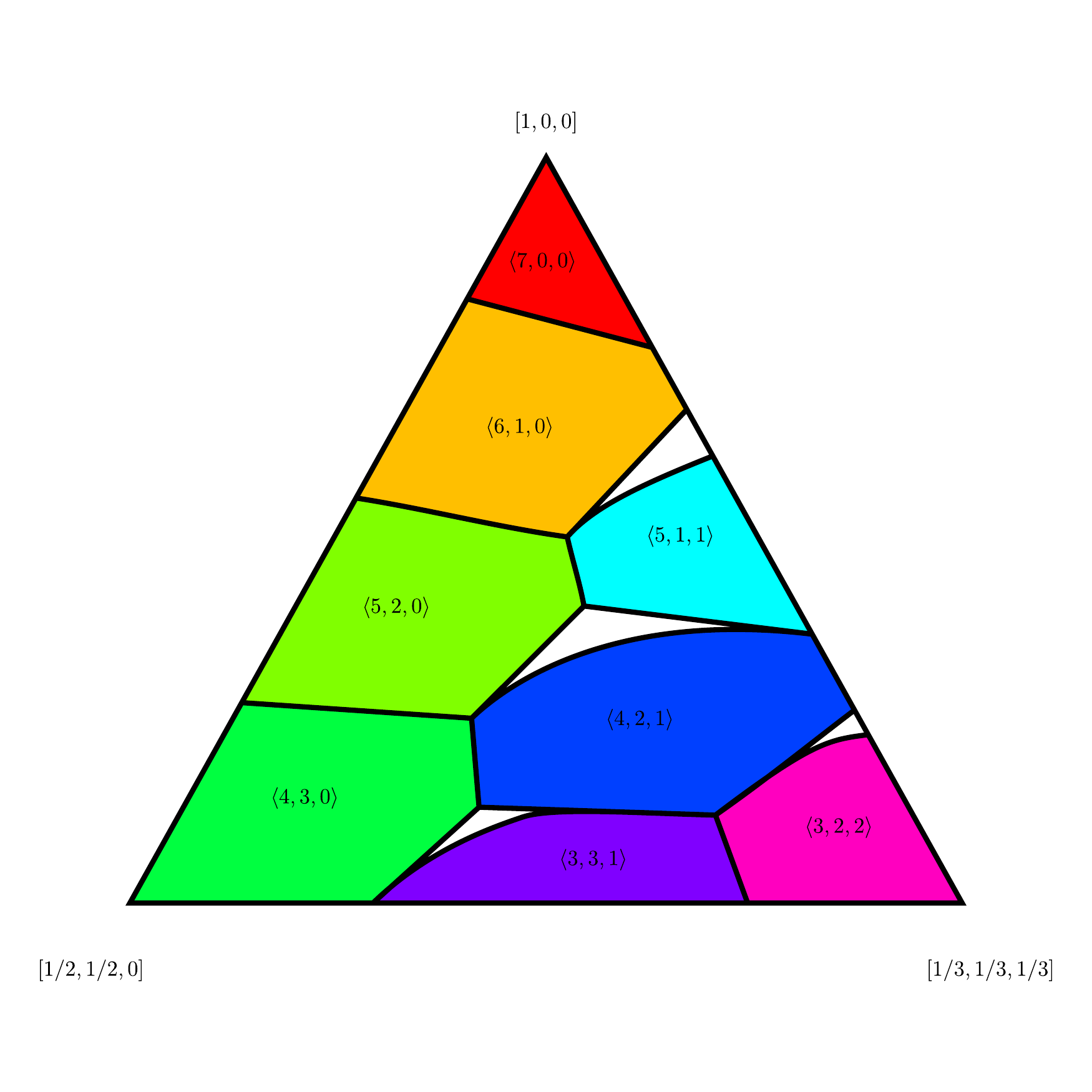}  
\caption{Triangular plot \label{winner_7} showing regions in which the
  minimax strategy is pure, marked in eight different colours for the
  eight different strategies~$\left<7,0,0\right>$
  through~$\left<3,2,2\right>$.  White signifies regions in which the
  minimax strategy is mixed}
\end{figure}

Figure~\ref{winner_7} illustrates the optimal strategy in a
game-theoretic sense for the case of~7 counters and shows that, for
some probabilities, the minimax strategy is not pure.  Thus if, for
example, $\left(p_1,p_2,p_3\right) =
\left(\frac{3}{4},\frac{1}{8},\frac{1}{8}\right)$, a minimax
strategy is to play~$\left<7,0,0\right>$ with probability 0.156,
$\left<6,1,0\right>$ with probability 0.189, and $\left<5,1,1\right>$
with probability 0.655.



\section{The game with an arbitrary number of boxes}

The PMF for the exact number of throws required to remove the final
counter in the general case is considerably more complicated.  We have
that~$\Prob(X=r)$ is

\begin{equation}
\sum_{i=1}^k \sum_{\sum_{j}r_j=r\atop r_i=0}
  {n+r-1 \choose
    n_1+r_1,\ldots,n_{i-1}+r_{i-1},n_i-1,n_{i+1}+r_{i+1},\ldots,n_k+r_k}
  \prod_{s=1}^k  p_s^{n_s+r_s}
  \label{pmf_general_case}
\end{equation}

Here the~$r_i$ are the number of ``wasted'' throws in box~$i$.  The
outer summation is over the number of the box whose counter is removed
last; the inner summation is over all distributions of~$r$ wasted
throws among the other boxes.

The original case~\cite{hankin2014} was to place~11 counters on boxes
corresponding to the total of two independent six-sided dice, that is,
with probabilities~$\min\left[i,13-i\right]/36,i=2,\ldots,12$.  There
are~56 ways in which the counters may be arranged if $p_i\geqslant
p_j\longrightarrow a_i\geqslant a_j$ is required.  Of these,
$\left<0,0,1,1,2,3,2,1,1,0,0\right>$ is optimal in expectation and
game-theoretic value.  There does not seem to be any way of
ascertaining this fact other than direct numerical evaluation of
all~56 expectations.

A more interesting example would arise from Zipf's
law~\cite{zipf1935}.  If there are~4 boxes then probabilities
proportional to~$\frac{1}{1},\frac{1}{2},\frac{1}{3},\frac{1}{4}$
would be indicated; with seven counters, there are~11 distinct
strategies from~$\left<7,0,0,0\right>$ through~$\left<2,2,2,1\right>$.

Expected time to removal is minimized by~$\left<5,1,1,0\right>$, while
the separate throw game has minimax strategy
of~$\left<4,2,1,0\right>$.  With a common throw, the minimax strategy
is mixed; specifically, play~$\left<4,1,1,1\right>$ with
probability~$0.26$ and $\left<3,2,1,1\right>$ with probability~$0.73$.

\section{Conclusions and Further work}

The liability allocation problem described here bears the same
relation to Blotto-type resource allocation as chore division does to
cake division.  The games discussed here do not appear to have been
discussed in the literature, even though they are simple to state, and
potentially rich in applications.  They are unusual in that the
simplest non-trivial cases are algebraically involved; numerical
methods have to be used for even small numbers of counters and boxes.

Further work might include relaxing the assumption that the box
probabilities are common knowledge; perhaps one or both players have
to infer the box probabilities from previous calls, as in multi-armed
bandit problems~\cite{berry1985}.

The three-box case shows a variety of results using the triplot
device; it would be interesting to understand why the boundary lines
are so close to perfect straight lines.

\subsection*{Acknowledgements}
The author thanks S. Marshall for valuable discussions; and also
A. M. Hankin (``Alice'') for bringing this problem to his attention,
making the (then) unintuitive observation that placing all one's
counters on the maximal probability box is suboptimal in expectation
and game value.

\bibliographystyle{spmpsci}      
\bibliography{manuscript_second_revision}   

\newpage

\begin{table}[h]
\newcolumntype{d}[1]{D{.}{\cdot}{#1} }
\begin{tabular}{c | d{2} d{2} d{2} d{2} d{2} d{2} }
I~$\backslash$~II &
\multicolumn{1}{r}{\strat{0}{5}} &
\multicolumn{1}{r}{\strat{1}{4}} &
\multicolumn{1}{r}{\strat{2}{3}} &
\multicolumn{1}{r}{\strat{3}{2}} &
\multicolumn{1}{r}{\strat{4}{1}} &
\multicolumn{1}{r}{\strat{5}{0}}  \\
\hline
\strat{0}{5}  & 0.00 & -0.84 & -0.53 & -0.16  &  0.19 & 0.47\\
\strat{1}{4}  & 0.84 &  0.00 & -0.33 &  0.09  &  0.42 & 0.65\\
\strat{2}{3}  & 0.53 &  0.33 &  0.00 &  0.37  &  0.64 & 0.81\\
\strat{3}{2}  & 0.16 & -0.09 & -0.37 &  0.00  &  0.83 & 0.92\\
\strat{4}{1}  &-0.19 & -0.42 & -0.64 & -0.83  &  0.00 & 0.98\\  
\strat{5}{0}  &-0.47 & -0.65 & -0.81 & -0.92  & -0.98 & 0.00\\
\end{tabular}
\caption{Expected game value for the two-box,
  five-counter \label{payoff_matrix_5counters_p=0.4} game, $p=0.4$,
  common throw.  Saddle point at \strat{2}{3}.}
\end{table}

\newpage

\begin{table}[h]
\begin{tabular}{c | ccc}
I~$\backslash$~II &
$\left<300\right>$&
$\left<210\right>$&
$\left<111\right>$\\
\hline
$\left<300\right>$ & 0 &
$\frac{2p_1^3}{(1-p_3)^3}-1$ &
$\begin{array}{l}2p_1^3\left(\frac{1}{(1-p_2)^3}\right.\\ \left. + \frac{1}{(1-p_3)^3}-1\right)-1\end{array}$
\\
$\left<210\right>$ & 
$1 - \frac{2p_1^3}{(1-p_3)^3}$ & 0 &\rule{0mm}{10mm}
$\begin{array}{l}p_2^2\left(\frac{2}{(1-p_2)^2}+\frac{1}{(1-p_3)^2}\right.\\ \qquad\left. +\frac{p_2}{1-p_1}-2\right)-1\end{array}$\\
$\left<111\right>$ & 
$\begin{array}{l}
1-2p_1^3\left(\frac{1}{(1-p_2)^3}\right.\\ \left. + \frac{1}{(1-p_3)^3}-1\right)\end{array}$
&
$\begin{array}{l}
1-p_2^2\left(\frac{2}{(1-p_2)^2}\right.\\ \left.+\frac{1}{(1-p_3)^2}+\frac{p_2}{1-p_1}-2\right)\end{array}$
&0
\\
\end{tabular}
\caption{Expected game value for the three-box, three counter common
  throw \label{payoff_matrix_3_boxes_3counters_common} game with
  probabilities~$\left(p_1,p_2,p_3\right)$}
\end{table}

\clearpage

\begin{table}[h]
\begin{tabular}{c | ccc}
I~$\backslash$~II &
$\left<300\right>$&
$\left<210\right>$&
$\left<111\right>$\\
\hline
$\left<300\right>$ & 0 & \rule{0mm}{6mm}
-0.059
&
+0.595\\
$\left<210\right>$ & \rule{0mm}{6mm}
+0.059
& 0 & +0.483\\
$\left<111\right>$ & \rule{0mm}{6mm}
-0.595 & -0.483& 0\\
\end{tabular}\hfill
\begin{tabular}{c | ccc}
I~$\backslash$~II &
$\left<300\right>$&
$\left<210\right>$&
$\left<111\right>$\\
\hline
$\left<300\right>$ & 0 & \rule{0mm}{6mm} +0.200 & +0.707\\
$\left<210\right>$ & \rule{0mm}{6mm} -0.200 & 0 & +0.519\\
$\left<111\right>$ & -0.707 &\rule{0mm}{6mm} -0.519 & 0\\
\end{tabular}
\caption{Expected game value, to three decimal places, for the
  three-box, three counter
 \label{numeric_payoff_matrix_3counters_3_boxes_p=0.6_0.3_0.1} game with
probabilities~$\left(\frac{6}{10},\frac{3}{10},\frac{1}{10}\right)$.
Left, common throw game; right, separate throw game}
\end{table}

\clearpage

\begin{table}[h]
\newcolumntype{d}[1]{D{.}{\cdot}{#1} }
\begin{tabular}{c | d{4 }}
strategy & \mathbb{E}X\\ \hline
$\left< 10,  0,  0\right>$ &   4.2840\\
$\left<  9,  1,  0\right>$ &   3.3731\\
$\left<  8,  2,  0\right>$ &   3.9637\\
$\left<  8,  1,  1\right>$ &   5.4261\\
$\left<  7,  3,  0\right>$ &   6.3922\\
$\left<  7,  2,  1\right>$ &   6.0878\\
$\left<  6,  4,  0\right>$ &  10.4550\\
$\left<  6,  3,  1\right>$ &   8.6927\\
$\left<  6,  2,  2\right>$ &  12.2692\\
$\left<  5,  5,  0\right>$ &  15.0561\\
$\left<  5,  4,  1\right>$ &  12.2627\\
$\left<  5,  3,  2\right>$ &  13.7832\\
$\left<  4,  4,  2\right>$ &  16.6344\\
$\left<  4,  3,  3\right>$ &  21.8522\\
\end{tabular}  
\caption{Expected time to removal for the three box, ten-counter game with
  probabilities~$(0.7,0.2,0.1)$ as calculated numerically;
 \label{expected_value_3boxes_10counters_table}
figures accurate to 4 decimal places}
\end{table}

\clearpage


\end{document}